\documentclass[12pt]{article}
\usepackage{amsmath}
\usepackage{amssymb}
\usepackage{cite}
\newsymbol \blackbox 1004
\newcommand{\eh}{\hfill}\newlength{\sperr}

\newenvironment{proof}{{\settowidth{\sperr}{\bf\rm
Proof}%
\par\addvspace{0.3cm}\noindent\parbox[t]{1.3\sperr}
{\bf\rm P\eh r\eh o\eh o\eh f\eh }%
}}{\nopagebreak\mbox{}
$\blackbox$\par\addvspace{0.3cm}}

\def\nn{\nonumber}
\def\a{\alpha}

\def\d{\delta}
\def\De{\Delta}

\def\vk{\varkappa}

\def\s{\sigma}
\def\la{\lambda}

\def\Om{\Omega}

\def\t{\theta}

\def\wh{\widehat}
\def\wt{\widetilde}
\def\ov{\overline}

\def\p{\partial}

\def\BC{{\mathbb C}}
\def\BR{{\mathbb R}}
\def\BN{{\mathbb N}}

\def\clq{{\mathcal Q}}

\def\cls{{\mathcal S}}

\def\diag{\mathrm{diag}}
\def\rk{\mathrm{Rank}}
\def\col{\mathrm{col}}
\newcommand{\E}{\mathrm{e}}
\newcommand{\I}{\mathrm{i}}

\newtheorem{Pa}{Paper}[section]
\newtheorem{Tm}[Pa]{{\bf Theorem}}

\newtheorem{Cy}[Pa]{{\bf Corollary}}
\newtheorem{Rk}[Pa]{{\bf Remark}}

\newtheorem{Ee}[Pa]{{\bf Example}}
\newtheorem{Dn}[Pa]{{\bf Definition}}
\newtheorem{Pn}[Pa]{{\bf Proposition}}

\title{Time-dependent Schr\"odinger equation in dimension $k+1$:
explicit and rational solutions via GBDT and multinodes}

\author{ A.L. Sakhnovich}

\date{}
\begin{document}
\maketitle

\thanks{E-mail address: Oleksandr.Sakhnovych@univie.ac.at}

\begin{abstract}  A version of the binary Darboux transformation is constructed
for  non-stationary Schr\"odinger equation  in dimension $k+1$, where $k$ is the number  of space variables,  $k \geq 1$. 
This is an iterated GBDT version. New families of non-singular
and rational potentials and solutions are obtained. Some results are new
for the case that $k=1$ too.
A certain generalization of a colligation introduced by M.S. Liv\v{s}ic and
of the $S$-node introduced by L.A. Sakhnovich is successfully used in our construction.

\end{abstract}


\vspace{1em}


\section{Introduction} \label{intro}
\setcounter{equation}{0}
The time-dependent, or non-stationary, Schr\"odinger (TDS) equation
is one of the most 
well-known and important equations in physics. We shall consider
the vector  case of the  Schr\"odinger equation
and an arbitrary coefficient $\a$
($\a \in \BC$, $\, \a \not=0$):
\begin{align}&      \label{1.1}
Hu=0, \quad H:=\a \frac{\p}{\p t}+\De- q(x,t), \quad u\in \BC^p, \end{align}
where  $\De$ denotes the Laplacian with respect to the spatial variables $x=(x_1,x_2, \ldots , x_k)\in
\BR^k$, $q$ is a $p \times p$ matrix function, $\BC$ stands for the complex plane, and $\BR$ stands for 
the real axis.
Explicit solutions of \eqref{1.1} are of essential and permanent
interest and various methods were applied to construct them.
For example, rational slowly decaying  soliton solutions (lumps) of the 
Kadomtsev-Petviashvili   (KP) 
and time-dependent  Schr\"odinger (TDS) equations were first found in \cite{MZB}. Since then they were actively studied as well
as rational solutions of other integrable equations (see, e.g., 
\cite{ACTV, EPV, FrN, Kri, ALSZ, VA, VPE, ZuS} and references therein).

Following seminal works of B\"acklund, Darboux, and Jacobi, different kinds of Darboux transformations and related commutation
and factorization methods are fruitfully used to obtain explicit solutions
of linear and nonlinear equations (see, e.g., \cite{Ci, cr, de, fg, GeH, gt, Gu, KoSaTe,
kr57, Mar, MS, Ros, ALS10, ZM} and numerous references therein).
In particular, matrix and operator identities were widely used in these constructions (see \cite{BGK, CSSC, DimM, FKS0, KaS, Mar, ALS94, ALS01,
ALS03, ALS06, ALSZ, Sch} for various results, discussions and references).

Darboux transformations proved to be especially useful  for
the construction of explicit solutions of the TDS
equation  in dimension $1+1$
(and, correspondingly, for the construction of explicit solutions of 
the KP  equation with two spatial variables).
A singular (non-binary) Darboux transformation was used for that purpose in
\cite{Mat} and a  binary Darboux transformation
for the scalar TDS appeared in the well-known book \cite[Section 2.4]{MS}.
Further important results on the
TDS  equation  in dimension $1+1$
can be found in \cite{ACTV, BaSa, BPPP, VA} (see also \cite{S-HPS}
and references therein for generalized TDS equations).

The first discussion on the Darboux transformation for TDS with $k>1$
spatial variables that we could find, was in the work \cite{Sab0} (Sabatier, 1991).
In spite of interesting publications (e.g, \cite{ABI, Sab}), the case
of linear equations with $k$ spatial variables ($k>1$) is much more difficult
(and much less is done, even for the singular Darboux transformation),
especially so for $k>2$ (see \cite{Sab, S-H} for some explanations).

In this paper we construct explicit and rational solutions of the TDS equation
for $k>1$ and $k>2$ as well. We apply the generalized B\"acklund-Darboux approach
(GBDT) from \cite{FKS0, KaS, ALS94, ALS01,
ALS03, ALS06, ALSZ} (see further references and some comparative
discussions on this method in \cite{Ci, KoSaTe, ALS10}).
Corresponding results for the case $k=1$ were announced in \cite{ALS03}
and the present paper contains proofs, which are valid for $k=1$ too.
GBDT is partially based on the operator identity (also called $S$-node) method
\cite{SaL1, SaL2, SaL3}, which, in its turn, takes its roots in the
characteristic matrix function and operator colligations introduced
by M.S. Liv\v{s}ic \cite{Liv1} (see also \cite{Liv2}).  In his later works
M.S. Liv\v{s}ic studied a greatly more complicated case
of colligations with several (instead of one) commuting 
operators \cite{Liv3} (see also \cite{BB, BaV, LKMV} and references therein).
Correspondingly, for the case of $k$ spatial variables we need an
$S$-node with $k$ matrix identities, which we call $S$-{\it multinode}.

In Section \ref{GBDT} we describe GBDT for the  TDS equation.
$S$-multinodes are introduced in Section \ref{ExplS}.
We use GBDT and multinodes in Section \ref{ExplS} to construct
explicitly solutions and potentials of the TDS equation and consider
examples. Conditions for non-singular
and rational solutions and potentials
and concrete examples are given in Section \ref{Lump}.

We use $\BN$ to denote the set of natural numbers, $\s$ to denote
spectrum,  and $\ov \a$ to denote the complex conjugate of $\a$.
The notatation $\rk(A)$  stands for the rank of a matrix $A$, $A^*$ is the
matrix adjoint to $A$,  and $\col$ denotes a column.

\section{GBDT for the TDS equation} \label{GBDT}
\setcounter{equation}{0}
Let $H=\a \frac{\p}{\p t}+\De- q$ be some TDS equation, which we call {\it initial},
and let $\Psi(x,t)$ and $\Phi(x,t)$ be  block rows of $n \times p$ blocks
$\Psi_r$ and $\Phi_r$, respectively. Here $n \in \BN$ is fixed
and $0\leq r \leq k$. It is required
that $\Psi$ satisfies equations
\begin{align}\label{2.1}&
H \Psi_0^*=0, \quad \Psi_r^* =\frac{\p}{\p x_r}\Psi_0^* \quad (1 \leq r \leq k),
\end{align}
where $H$ is applied to $ \Psi_0^*$ columnwise. In other words, $\Psi^*$ satisfies a first order differential system,
which is equivalent to TDS:
\begin{align}\label{2.2}&
L \Psi^*=0, \quad L:=\a \t_0 \frac{\p}{\p t}+\sum_{r=1}^k\t_r \frac{\p}{\p x_r} -\t_{k+1};
\\ \nn &
\t_0:=e_k e_0^*, \qquad  \t_r:=e_{r-1}e_0^* +e_ke_r^* \quad (1\leq r\leq k), 
\\ \label{2.3}&
\t_{k+1}(x,t):=e_k q(x,t) e_0^*+ \{\d_{i,j-1}I_p\}_{i,j=0}^k,
\end{align}
where $\, \d_{i,j}\,$ is the Kronecker's delta, $I_p$ is the $p \times p\,$ identity
matrix, and $\, (k+1)p \times p$ matrices $e_i$ are given by the 
equilities $e_i:= \{\d_{i,j}I_p\}_{j=0}^k$. We require
\begin{align}\label{2.1'}&
 \Phi_r =\frac{\p}{\p x_r}\Phi_0 \qquad (1 \leq r \leq k),
\end{align}
and $\Phi_0$ shall  be discussed a bit later.

An $n \times n$ matrix function ${\cal S}$, which we define via $\Psi$ and $\Phi$:
\begin{align}\label{2.4}&
\frac{\p}{\p x_r}{\cal S}(x,t)=\Phi_0(x,t) \nu_r \Psi_0(x,t)^* \qquad (1\leq r\leq k), \\
\label{2.5}&
\frac{\p}{\p t}{\cal S}(x,t)=
\a^{-1}\sum_{r=1}^k\big(\Phi_r(x,t)\nu_r\Psi_0(x,t)^*-\Phi_0(x,t)\nu_r\Psi_r(x,t)^*\big),
\end{align}
 is very important in GBDT.  Here $\nu_r$ are some $p\times p$ matrices. 
For linear equations depending on one variable
and nonlinear equations depending on two variables, the analog 
of ${\cal S}$ is denoted
by $S$ and the so called Darboux matrix is presented as the
transfer matrix function of the corresponding $S$-node. The equality
$\wt \Psi_0^*= \Psi_0^*{\cal S}^{-1}$
for the solution  $\wt \Psi_0^*$ of the transformed TDS $\wt Hf=0$
holds also in our case (see Theorem \ref{TmGBDT} below).

Since GBDT is a kind of binary Darboux transform, the matrix function $\Phi_0$
should satisfy some dual to TDS differential equation $H_d\Phi_0 =0$. 
In view of \eqref{2.1}, \eqref{2.1'}, \eqref{2.4}, and \eqref{2.5}, for the case
that
\begin{align}\label{2.7}&
H_df=\a \frac{\p}{\p t}f(x,t)-\De f(x,t)- f(x,t) q_d(x,t), \qquad q_d\nu_i =\nu_i q \quad (1\leq i\leq k),
\end{align}
where $f$ is a row vector function and $q_d$ is a $p\times p$ matrix function, we have
\begin{equation}\label{2.8}
{\cal S}_{tx_i}=\a^{-1}\sum_{r=1}^k\Big(\Phi_r \nu_r \Psi_i^*  - \Phi_i \nu_r \Psi_r^*+
\Big(\frac{\p}{\p x_i}\Phi_r\Big)\nu_r\Psi_0^*-\Phi_0\nu_r\Big(\frac{\p}{\p x_i}\Psi_r^*\Big)\Big),
\end{equation}
\begin{align}
\nn
{\cal S}_{x_i t}&=
\a^{-1}\big((\De \Phi_0)\nu_i\Psi_0^*-\Phi_0\nu_i(\De \Psi_0^*)+\Phi_0(\nu_i q -q_d\nu_i) \Psi_0^*\big)
\\
\label{2.9}&
=\a^{-1}\big((\De \Phi_0)\nu_i\Psi_0^*-\Phi_0\nu_i(\De \Psi_0^*)\big), \qquad {\cal S}_{x_i t}:=\frac{\p}{\p t} \frac{\p}{\p x_i}{\cal S}.
\end{align}
Because of \eqref{2.8} and \eqref{2.9}, the compatibility condition ${\cal S}_{x_1 t}={\cal S}_{t x_1}$ for equations
\eqref{2.4} and \eqref{2.5}  is fulfilled
for the case that $k=1$. However, for $k>1$ the situation is more complicated, and 
we just assume the existance of ${\cal S}$, satisfying \eqref{2.4} and \eqref{2.5}, and don't assume \eqref{2.7} in our theorem below.
\begin{Tm} \label{TmGBDT} Let matrix functions $\Psi$, $\Phi$, and ${\cal S}$
satisfy relations \eqref{2.1}, \eqref{2.1'}, and \eqref{2.4}, \eqref{2.5},
respectively. Then, in the points of invertibility of ${\cal S}$, the matrix function
\begin{align}\label{2.6}&
\wt \Psi_0^*:= \Psi_0^*{\cal S}^{-1}
\end{align}
 satisfies the transformed TDS equation:
 \begin{align}\label{2.11}&
\wt H\wt \Psi_0^*=0, \quad \wt H:=\a \frac{\p}{\p t}+\De- \wt q(x,t),
\end{align}
where 
 \begin{align}\label{2.12}&
 \wt q(x,t):=q(x,t)-2\sum_{r=1}^k 
 \frac{\p}{\p x_r}\big(\Psi_0(x,t)^*{\cal S}(x,t)^{-1}\Phi_0(x,t)\big)\nu_r.
\end{align}
\end{Tm}
\begin{proof}. Taking into account \eqref{2.1} and definitions
of $H$, $\Psi_0$, and $\wt H$ in \eqref{1.1}, \eqref{2.6}, and \eqref{2.11},
respectively, we get
 \begin{align}\label{2.13}&
\wt H\wt \Psi_0^*=(q-\wt q)\wt \Psi_0^*- \a \Psi_0^*\cls^{-1}\cls_t \cls^{-1}+
\Psi_0^*\De(\cls^{-1})- 2 \sum_{r=1}^k\Psi_r^* \cls^{-1}\cls_{x_r} \cls^{-1}.
\end{align}
Because of  \eqref{2.4}, we have
 \begin{align}\label{2.14}&
\De(\cls^{-1})= \sum_{r=1}^k \cls^{-1}\big(
2\cls_{x_r} \cls^{-1}\cls_{x_r} -\Phi_r\nu_r\Psi_0^*-\Phi_0\nu_r\Psi_r^*
\big)\cls^{-1}.
\end{align}
Finally, using formulas \eqref{2.5} and \eqref{2.14} and reducing
similar terms, we rewrite  \eqref{2.13} as
 \begin{align}\nn
\wt H\wt \Psi_0^*=&(q-\wt q)\wt \Psi_0^*-\Psi_0^*\cls^{-1}
\sum_{r=1}^k\Phi_r\nu_r
\wt \Psi_0^*+2\Psi_0^*\cls^{-1}\Phi_0 \sum_{r=1}^k\nu_r\Psi_0^*\cls^{-1}\Phi_0
\nu_r \wt \Psi_0^*
\\ \label{2.15}&
-\Psi_0^*\cls^{-1}\sum_{r=1}^k\Phi_r\nu_r\wt \Psi_0^*
-2\sum_{r=1}^k\Psi_r^*\cls^{-1}\Phi_0\nu_r\wt \Psi_0^*.
\end{align}
Since $\Phi_r=(\Phi_0)_{x_r}$, $\Psi_r=(\Psi_0)_{x_r}$, and
\eqref{2.4} holds,
it follows from \eqref{2.15} that for $\wt q$ given by \eqref{2.12}
the equality $\wt H\wt \Psi_0^*=0$ is true.
\end{proof}

\section{Multinodes and explicit solutions} \label{ExplS}
\setcounter{equation}{0}
\begin{Dn}\label{DnS}
By a matrix $S$-multinode $($or, more precisely, by $S_k$-node
$\big\{k,A,B,R,\nu,C_{\Phi},C_{\Psi}\big\}$$)$
we call a set of matrices, which consists of $N \times N$ commuting
matrices $A_r\,$ $( 1\leq r\leq k)$,
of  $N \times N$ commuting  matrices $B_r\,$ $( 1\leq r\leq k)$,
of $p \times p$ matrices $\nu_r\,$ $( 1\leq r\leq k)$, and of
an $N \times N$ matrix $R$, an $N \times p$ matrix $C_{\Phi}$,
and a $p \times N$ matrix $C_{\Psi}$, such that the matrix identities
\begin{align}\label{3.1}&
A_rR-RB_r=C_{\Phi}\nu_rC_{\Psi}, \qquad 1\leq r\leq k
\end{align}
hold. An operator $S$-multinode is defined in the same way.
\end{Dn}

For the case that $k=1$ this definition coincides with the definition
of an $S$-node \cite{SaL1, SaL2, SaL3}, and for the case that $R=I_N$ 
and $B_r=A_r^*$ our definition coincides with the definition
of a colligation from \cite{Liv3}.

In this section we treat the case $q\equiv 0$, that is,
\begin{align}\label{3.1'}&
 \wt q(x,t)=-2\sum_{r=1}^k 
 \frac{\p}{\p x_r}\big(\Psi_0(x,t)^*{\cal S}(x,t)^{-1}\Phi_0(x,t)\big)\nu_r.
\end{align}
\begin{Tm}\label{TmExplS} Let an $n \times N$ matrix $\wh C_{\Phi}$,
an $N \times n$ matrix $\wh C_{\Psi}$, an $n \times n$ matrix $\cls_0$,
and a matrix $S_k$-node 
$\big\{k,A,B,R,\nu,C_{\Phi},C_{\Psi}\big\}$ be given.
Then the matrix functions 
\begin{align}\label{3.2}&
\Phi_0(x,t)=\wh C_{\Phi}\E_A(x,t)C_{\Phi}, \quad 
\E_A(x,t):=\exp \Big\{\big(\sum_{r=1}^k x_rA_r\big)+
\a^{-1}t\big(\sum_{r=1}^k A_r^2\big)\Big\},
\\ \label{3.3}&
 \Psi_0(x,t)^*= C_{\Psi}\E_B(-x,-t)\wh C_{\Psi}, \quad
{\cal S}=\wh C_{\Phi}\E_A(x,t)R\E_B(-x,-t)\wh C_{\Psi}+{\cal S}_0
\end{align}
satisfy conditions of Theorem \ref{TmGBDT}, where the initial TDS equation
is chosen so that $q\equiv 0$ $($i.e., the conditions on $\Phi$ and $\Psi$
are satisfied after we standardly add $\Phi_r=(\Phi_0)_{x_r}$ and
 $\Psi_r=(\Psi_0)_{x_r}$$)$.
\end{Tm}
\begin{proof}. It is immediate from \eqref{3.3} that $\a \frac{\p}{\p t}\Psi_0^*+\De \Psi_0^*=0$, that is,
$H\Psi_0^*=0$, and so \eqref{2.1} holds. Because of \eqref{3.1}, \eqref{3.2}, and
\eqref{3.3} we have \eqref{2.4}.
 It remains to show that \eqref{2.5} holds. For that purpose, note that equalities \eqref{3.1} imply
\begin{align}\nn 
A_r^2R-RB_r^2 & =A_r(A_rR-RB_r)+(A_rR-RB_r)B_r
\\ \label{3.4}&
=A_rC_{\Phi}\nu_rC_{\Psi}+C_{\Phi}\nu_rC_{\Psi}B_r, \quad   1\leq r\leq k.
\end{align}
In view of \eqref{3.3} and \eqref{3.4} we get
\begin{align}\label{3.5}&
\frac{\p}{\p t}\cls=\a^{-1}\wh C_{\Phi}\E_A(x,t)\Big(\sum_{r=1}^k 
A_rC_{\Phi}\nu_rC_{\Psi}+C_{\Phi}\nu_rC_{\Psi}B_r\Big) \E_B(-x,-t)\wh C_{\Psi}.
\end{align}
Since $\Phi_r=(\Phi_0)_{x_r}$ and
 $\Psi_r=(\Psi_0)_{x_r}$, formula \eqref{2.5} easily follows from \eqref{3.2}, \eqref{3.3}, and \eqref{3.5}.
\end{proof}
\begin{Rk}\label{RkPhi} 
Clearly, the equality $\a \frac{\p}{\p t}\Phi_0- \De \Phi_0=0$ holds for
$\Phi_0$ of the form \eqref{3.2},
that is, we have $H_d\Phi_0=0$
for the case that $q_d=0$, and so \eqref{2.7} holds. Hence, it goes
in Theorem \ref{TmExplS} about the binary GBDT.
\end{Rk}
Recall that singular (and some stationary binary) Darboux transformations
of a scalar TDS equation into the vector TDS were treated in \cite{Sab}.
It is easy to see that for the case that 
$\wt q=\{\wt q_{ij}\}_{i,j=1}^p$ and  $f=\{ f_{i}\}_{i=1}^p$ are the potential
and solution, respectively, of some  vector TDS equation, 
the functions 
\begin{align}\label{sc}&
q^{sc}=\wt q_{rr}+\big(\sum_{j\not=r}\wt q_{rj}f_j\big)/f_r, \quad
f^{sc}=f_r
\end{align}
are the potential
and solution of a scalar TDS.
\begin{Ee}\label{Eek1} Consider the simplest case $k=1$ $($and $\nu_1=I_p)$.
If $A_1$ and $B_1$ are diagonal matrices:
$A_1=\diag\{a_1,\ldots,a_N\}$ and $B_1=\diag\{b_1,\ldots,b_N\}$, then identity
\eqref{3.1} leads us to 
\begin{align}\label{3.6}&
R_{ij}=(a_i-b_j)^{-1}\big(C_{\Phi}C_{\Psi}\big)_{ij}, \quad
R=\{R_{ij}\}_{i,j=1}^N,
\end{align}
and relations \eqref{2.6}, \eqref{3.1'}-\eqref{3.3}, and \eqref{3.6}
define  solution and potential
of the transformed vector TDS equation explicitly $($up to inversion of matrices$)$.
\end{Ee}
\begin{Ee}\label{Eekm1}  Now, consider a simple example for the case 
that $1<k \leq p$.  We set
\begin{align}\label{3.7}&
A_1=\diag\{a_1,\ldots,a_N\}, \quad A_r=(c_r I_N-A_1)^{-1} \quad (a_i \not=c_r,
\, \,1<r \leq k);
\\  \label{3.8}&
B_1=\diag\{b_1,\ldots,b_N\}, \quad B_r=(c_r I_N-B_1)^{-1} \quad 
(b_i \not=c_r,
\, \,1<r \leq k);
\\ \label{3.9}&
R=\{(a_i-b_j)^{-1}\}_{i,j=1}^N \,\, (a_i \not= b_j); \quad \nu_r=\{\d_{r-i}\d_{r-j}\}_{i,j=1}^p, 
\,\, 1\leq r \leq k;
\\ \label{3.10}&
C_{\Phi}=\begin{bmatrix}
h & A_2 h & \ldots & A_kh&\breve C_{\Phi}
\end{bmatrix}, \quad h=\col \, \begin{bmatrix}
1 & 1  & \ldots & 1
\end{bmatrix},
\\ \label{3.11}&
C_{\Psi}=\col \, \begin{bmatrix}
h^* & h^* B_2  & \ldots & h^*B_k& \breve C_{\Psi}
\end{bmatrix},
\end{align}
where $\col$ stands for column, $\breve C_{\Phi}$ is some $N \times (p-k)$
matrix, and $\breve C_{\Psi}$ is some $(p-k) \times N$
matrix.  We see that $A_1R-RB_1=hh^*$. Hence, the equalities
$A_rR-RB_r=A_r h h^*B_r$ $(1<r \leq k)$ hold. Therefore,
relations  \eqref{3.7}-\eqref{3.11}  determine an
$S_k$-node $\big\{k,A,B,R,\nu,C_{\Phi},C_{\Psi}\big\}$ and 
the corresponding explicit expressions for solution and potential of the transformed TDS equation follow.
\end{Ee}
GBDT with non-diagonizable matrices $A_r$ is of interest (and has its own specifics). See, for instance, \cite{ALS03} where the cases that $k=1$
and matrices $A_1$ are $2 \times 2$ or $3\times 3$ Jordan cells are
treated as examples. Below we present an example of a multinode,
where $k>1$ and matrices $A_r$  are non-diagonizable.
\begin{Ee}\label{EeND} We assume $1 \leq k \leq p$, $N>1$ and set
\begin{align}\label{4.1}&
A_0:=\frac{\I}{2}I_N+\{a_{i-j}\}_{i,j=1}^N, \quad a_s:=\I \quad {\mathrm{for}}
\,\, s>0, \,\, a_s:=0 \quad {\mathrm{for}}
\,\, s\leq 0;
\\ \label{4.2}&
A_r:=(c_r I_N-A_0)^{-1}, \,\, B_r:=(c_r I_N-A_0^*)^{-1} \,\, (c_r \not=\pm \I /2,
\, \,1\leq r \leq k).
\end{align}
We see that matrix $A_0$, and therefore matrices $A_r$ and $B_r$ are linear similar to Jordan cells.
The matrix $R$ is a so called cyclic Toeplitz matrix and is introduced by the equality
\begin{align}\label{4.3}&
R:=\{T_{i-j}\}_{i,j=1}^N \quad (T_s \in \BC, \quad T_s=0 \quad {\mathrm{for}}
\quad s<0).
\end{align}
Then, the following matrix identity holds (see, e.g., \cite[p. 451]{ALS00}):
\begin{align}\label{4.4}&
A_0R-RA_0^*=igh^*, \quad g= \col \, \begin{bmatrix}
T_0 & T_0+T_1  & \ldots & \sum_{s=0}^{N-1}T_s
\end{bmatrix}, 
\end{align}
where $h$ is given in \eqref{3.10}. Because of \eqref{4.2} and \eqref{4.4},
the identities \eqref{3.1} are true, where $\nu_r$ are given in \eqref{3.9} and
\begin{align} \label{4.5}&
C_{\Phi}=\I \begin{bmatrix}
A_1g  & \ldots & A_k g&\breve C_{\Phi}
\end{bmatrix}, \quad 
C_{\Psi}=  \col \, \begin{bmatrix}
h^*B_1  & \ldots & h^*B_k& \breve C_{\Psi}
\end{bmatrix}.
\end{align}
That is,  a multinode,
where $k>1$ and matrices $A_r\,$  and $\, B_r$ are non-diagonizable,
is constructed.
\end{Ee}
The cases, where matrices $\nu_r$ had rank 1 and $p \times p$ matrix 
TDS equations with $p$ spatial variables were included, were treated
in Examples \ref{Eekm1} and \ref{EeND}. 
\begin{Rk}\label{rank}
Clearly, it is quite possible, though
somewhat less convenient, to consider $S_k$-nodes  with matrices $\nu_r$
of higher ranks in the same way. Recall also an easy transfer \eqref{sc}
from a matrix to a scalar TDS.
\end{Rk}
Definition \ref{DnS} admits an easy generalization for the case
of rectangular matrices $R$, whereupon the  proof
of Theorem \ref{TmExplS} does not require any changes.
\begin{Dn}\label{DnRect} 
By the $S_k$-node
$\big\{k,A,B,R,\nu,C_{\Phi},C_{\Psi}\big\}$ $($with rectangular matrix $R)$
we call a set of matrices, which consists of $N_1 \times N_1$ commuting
matrices $A_r\,$ $( 1\leq r\leq k)$,
of  $N_2 \times N_2$ commuting  matrices $B_r\,$ $( 1\leq r\leq k)$,
of $p \times p$ matrices $\nu_r\,$ $( 1\leq r\leq k)$, and of
an $N_1 \times N_2$ matrix $R$, an $N_1 \times p$ matrix $C_{\Phi}$,
and a $p \times N_2$ matrix $C_{\Psi}$, such that the matrix identities
\eqref{3.1} hold.
\end{Dn}
\begin{Cy}\label{CyExplS}
Let an $n \times N_1$ matrix $\wh C_{\Phi}$,
an $N_2 \times n$ matrix $\wh C_{\Psi}$, an $n \times n$ matrix $\cls_0$,
and a matrix $S_k$-node 
$\big\{k,A,B,R,\nu,C_{\Phi},C_{\Psi}\big\}$ $($with an $N_1 \times N_2$ matrix $R)$
be given.
Then the matrix functions $\Phi_0$, $\Psi_0$, and $\cls$,
which are given by formulas \eqref{3.2} and \eqref{3.3},
satisfy equation $H_d\Phi_0=0$ and conditions of Theorem \ref{TmGBDT},
where $q=q_d=0$.
\end{Cy}
\section{Non-singular, rational, and lump potentials} \label{Lump}
\setcounter{equation}{0}
In this section, we study conditions that the potentials $\wt q$ 
and the TDS solutions $\wt \Psi_0^*$
 are 
non-singular and rational. 
Our next proposition deals with a construction of rational potentials.
\begin{Pn}\label{PnRS}
Let an $n \times N_1$ matrix $\wh C_{\Phi}$,
an $N_2 \times n$ matrix $\wh C_{\Psi}$, an $n \times n$ matrix $\cls_0$,
and a matrix $S_k$-node 
$\big\{k,A,B,R,\nu,C_{\Phi},C_{\Psi}\big\}$ 
be given. Assume additionally that the conditions $(i)$ 
below hold:

$(i)$ all the matrices from the set 
$\{A_r\}\cup \{B_r\}$
are nilpotent.

Then the solution and potential
of the transformed TDS equation, which are given by formulas \eqref{2.6} and \eqref{3.1'}, respectively,
are rational.

If, instead of $(i)$, we assume that
\begin{align}\label{4.30}&
\cls_0=0; \quad A_r=\mu_rI_{N_1} +\breve A_r, \quad B_r=\la_rI_{N_2}+\breve B_r
,
\end{align}
where $1\leq r\leq k; \quad \mu_r, \, \la_r \in \BC$,
and the matrices $\breve A_r$ and $\breve B_r$ are nilpotent,
then the potential $\wt q$
of the transformed TDS equation is rational.
\end{Pn}
\begin{proof}.
In view of  Corollary \ref{CyExplS},  the matrix functions $\Phi_0$, 
$\Psi_0$, and $\cls$
satisfy  conditions of Theorem \ref{TmGBDT}. Therefore, using Theorem \ref{TmGBDT} we see that the solution and potential
of the transformed TDS equation are given by formulas \eqref{2.6} and \eqref{3.1'},
respectively.

First, assume that the conditions $(i)$ hold.
It is immediate that all the matrix functions 
\[\exp\{x_rA_r\}, \quad 
\exp\{\a^{-1}tA_r^2\}, \quad \exp\{-x_rB_r\}, \quad \exp\{-\a^{-1}tB_r^2\}
\quad (1\leq r \leq k)
\]
are matrix polynomials, and thus $\Phi_0$, $\Psi_0$, and $\cls$, which are given by \eqref{3.2} and \eqref{3.3},
are matrix polynomials with respect to $x$ and $t$.
The statement of the  proposition follows.

Next, assume that  conditions \eqref{4.30} hold. Because of  \eqref{4.30}
we have
\begin{align}\label{4.31}&
\E_A(x,t)=\E^{f(x,t)}p_A(x,t), \quad \E_B(-x,-t)=\E^{f(x,t)}p_B(x,t),
\end{align}
where $p_A$ and $p_B$ are matrix polynomials, whereas $f$ and $g$
are scalar polynomials:
\begin{align}\label{4.32}&
f(x,t)=\sum_{r=1}^k\mu_rx_r+\a^{-1}\sum_{r=1}^k\mu_r^2t,
\quad g(x,t)=-\sum_{r=1}^k\la_rx_r-\a^{-1}\sum_{r=1}^k\la_r^2t.
\end{align}
Since $\cls_0=0$, we derive from
formulas \eqref{3.2}, \eqref{3.3}, and \eqref{4.31} that $\Psi_0^*\cls^{-1}\Phi_0$
is rational, and so (in view of \eqref{3.1'}) the potential $\wt q$ is rational too.
\end{proof}
Further we assume again that $N_1=N_2=N$, that is, $R$ is a square matrix.
The following proposition is immediate from
\eqref{2.6}, \eqref{3.1'}-\eqref{3.3}.
\begin{Pn}\label{PnSA}
Let the conditions of Theorem \ref{TmExplS} hold and let also equalities
\begin{align}\label{4.6}&
\a=\I, \quad \wh C_{\Psi}= \wh C_{\Phi}^*, \quad \cls_0 = \cls_0^*;
\\ \label{4.7}&
C_{\Psi}= C_{\Phi}^*, \quad R=R^*; \quad B_r=-A_r^*
\quad (1\leq r \leq k) 
\end{align}
be satisfied. Then we have
\begin{align}\label{4.8}&
\E_B(-x, -t)=\E_A(x,t)^*, \quad \Phi_0(x,t)=\Psi_0(x,t), \quad \cls(x,t)=\cls(x,t)^*.
\end{align}
Furthermore, if the additional relations
\begin{align}\label{4.9}&
R \geq 0, \quad \cls_0>0 \quad {\mathrm{or}} \quad R > 0, \quad
\rk(\wh C_{\Phi})=n \leq N,
\quad \cls_0 \geq 0
\end{align}
hold, then  the inequality $\cls(x,t)>0$ holds too, and so $\cls$ is invertible
and the solution $\wt \Psi_0^*$ and potential $\wt q$ of the transformed TDS
are non-singular.
\end{Pn}
Finally, we consider several concrete examples of  non-singular, rational, and lump potentials, where \eqref{4.6} and  \eqref{4.7} hold, and
\begin{align}\label{4.10}&
\nu_r=\{\d_{r-i}\d_{r-j}\}_{i,j=1}^p, 
\qquad 1\leq r \leq k.
\end{align}

Because of \eqref{4.10} and the first and third
equalities in \eqref{4.7}, identity \eqref{3.1} for $r=1$ has the form
\begin{align}\nn &
A_1R+RA_1^*=\clq, \quad \clq=\clq^*, \quad {\mathrm{i.e.,}} 
\\ \label{4.11} &
(zI_N-\I A_1)^{-1}R
=
\I(zI_N-\I A_1)^{-1}\clq (zI_N+\I A_1^*)^{-1}+
R(zI_N+\I A_1^*)^{-1}.
\end{align}
If $\s(\I A) \subset \BC_+$, we take residues and derive from \eqref{4.11}
a well-known representation
\begin{align}\label{4.12}&
R=\frac{1}{2\pi}\int_{-\infty}^{\infty}(zI_N-\I A_1)^{-1}\clq (zI_N+\I A_1^*)^{-1}dz,
\quad R=R^*,
\end{align}
that is, the second equality in \eqref{4.7} follows now from the first and third
equalities.
\begin{Rk}\label{IR} For the case that $\clq \geq 0$, representation
\eqref{4.12} implies $R\geq 0$. 
\end{Rk}
\begin{Ee} \label{CE1}
Let $k=p=3$, $N=2$, and let $A_1$ be a $2 \times 2$ Jordan cell:
\begin{align}\label{4.13}&
A_0=\begin{bmatrix}
0 & 1 \\ 0 & 0 
\end{bmatrix}, \quad A_1=\mu_0 I_2+A_0, \quad \vk:= \mu_0+\ov{\mu_0}\not= 0;
\\ \label{4.14}&
A_r=(A_1-c_rI_2)^{-1},
\quad c_r=-\ov {c_r} \qquad (r=2,3).
\end{align}
It is immediate from \eqref{4.13} and \eqref{4.14} that
\begin{align}\label{4.15}&
A_r=(\mu_0-c_r)^{-1}I_2-(\mu_0-c_r)^{-2}A_0,
\end{align}
and $($similar to \eqref{3.10}  but with  different choices of $A_r$ and $h$$)$ we put 
\begin{align}\label{4.16}&
C_{\Psi}=C_{\Phi}^*, \quad C_{\Phi}=\begin{bmatrix}
h & A_2 h & A_3 h
\end{bmatrix}=\begin{bmatrix}
0 & -(\mu_0-c_2)^{-2} & -(\mu_0-c_3)^{-2}
\\
1 & (\mu_0-c_2)^{-1} & (\mu_0-c_3)^{-1}
\end{bmatrix},
\end{align}
where $h=\col \begin{bmatrix}
0 & 1  
\end{bmatrix}$. We recover $R$ from the identity
\begin{align}\label{4.17}&
A_1R+RA_1^*=hh^*,
\end{align}
which is equivalent $($in view of $B_1=-A_1^*$, \eqref{4.10}, and \eqref{4.16}$)$ to  relation \eqref{3.1} for $r=1$ .
That is, we rewrite \eqref{4.17} in the form
\begin{align}\label{4.18}&
\vk R+\begin{bmatrix}
R_{21}+R_{12} & R_{22}\\ R_{22} & 0  
\end{bmatrix}=hh^*,
\end{align}
where $R_{ij}$ are the entries of $R$. Using \eqref{4.18} $($and starting from recovery of $R_{22}$$)$,
we easily get a unique $R$ satisfying \eqref{4.17}$:$
\begin{align}\label{4.19}&
R=\begin{bmatrix}
2\vk^{-3} & -\vk^{-2}\\ -\vk^{-2} & \vk^{-1}  
\end{bmatrix}, \quad \det R=\vk^{-4}.
\end{align}
Identities \eqref{3.1} for $r>1$ easily follow from \eqref{4.18}, and so we obtain
an $S_3$-node $\{3, A,B, R, \nu, C_{\Phi},C_{\Psi}\}$. Moreover, Remark \ref{IR}
and the second equality in \eqref{4.19} yield
\begin{align}\label{4.20}&
R>0 \quad {\mathrm{for}} \quad \vk=\mu_0+\ov{\mu_0}>0.
\end{align}
Therefore, the conditions of Theorem \ref{TmExplS} and Proposition \ref{PnSA}
are fulfilled.

Since $A_0^2=0$, formula \eqref{4.15} and second relations in \eqref{3.2} and \eqref{4.13} imply
\begin{align}\label{4.21}
\E_A(x,t)=&\E^{\Om_0(x,t)}(I_2+\Om_1(x,t)A_0)=\E^{\Om_0(x,t)}\begin{bmatrix}
1 & \Om_1(x,t)\\ 0 & 1  
\end{bmatrix}, \\ \nn  \Om_0(x,t):=& \mu_0 x_1+ (\mu_0-c_2)^{-1}x_2+ (\mu_0-c_3)^{-1}x_3
\\  \label{4.22} &
-\I \big(\mu_0^2+
 (\mu_0-c_2)^{-2}+ (\mu_0-c_3)^{-2}\big)t, \\ \nn
  \Om_1(x,t):=& x_1- (\mu_0-c_2)^{-2}x_2
- (\mu_0-c_3)^{-2}x_3
\\ \label{4.23}&
-2\I \big(\mu_0-
 (\mu_0-c_2)^{-3}- (\mu_0-c_3)^{-3}\big)t.
\end{align}
\end{Ee}
\begin{Cy}\label{CyN2} Let $k=p=3$, $N=2$, and $ \vk=\mu_0+\ov{\mu_0}>0$.
Define matrices $C_{\Psi}$, $C_{\Phi}$, and $R$ via \eqref{4.16}
and \eqref{4.19}. Choose $\wh C_{\Psi}= \wh C_{\Phi}^*$ and $\cls_0>0$.
Then   relations \eqref{2.6} and \eqref{3.1'}-\eqref{3.3},
where $\E_A$ is given by \eqref{4.21}-\eqref{4.23} and $\E_B(-x, -t)=\E_A(x,t)^*$,
explicitly define non-singular solutions and potentials of TDS.
\end{Cy}
The cases $N>2$ can be treated in the same way.
\begin{Ee} \label{CE2} Let $k=p=3$ and $N=3$. Set 
\begin{align}\label{4.24}&
A_1=\mu_0 I_3+A_0, \quad A_r=(A_1-c_rI_3)^{-1}=
- \sum_{i=0}^2(c_r- \mu_0)^{-i-1}A_0^i,
\\ \label{4.25}& 
c_r=-\ov {c_r} \quad (r=2,3); \quad
 \vk:= \mu_0+\ov{\mu_0}\not= 0,  
\\ \label{4.26}& 
A_0=\begin{bmatrix}
0 & 1 & 0\\ 0 & 0 &1
\\ 0& 0& 0
\end{bmatrix}, 
\quad 
h=\begin{bmatrix}
0 \\ 0 \\ 1
\end{bmatrix}, \quad
C_{\Phi}=C_{\Psi}^*=\begin{bmatrix}
h & A_2 h & A_3 h
\end{bmatrix}.
\end{align}
In a way, which is quite similar to Example \ref{CE1}, 
we can show that an  $S_3$-node $\{3, A,B, R, \nu, C_{\Phi},C_{\Psi}\}$
appears, if we put
\begin{align}\label{4.27}&
B_r=-A_r^* \quad (1 \leq r \leq 3), \quad R=\begin{bmatrix}
6 \vk^{-5} & -3\vk^{-4} & \vk^{-3}
\\ -3 \vk^{-4} & 2\vk^{-3} & -\vk^{-2}\\ \vk^{-3} & -\vk^{-2} & \vk^{-1}  
\end{bmatrix}.
\end{align}
Moreover, we have $\det R=\vk^{-9}\not=0$, and so \eqref{4.20}
is valid for $R$ of the form \eqref{4.27} too.
Finally, we note that since $A_0^3=0$, formulas  \eqref{3.2} and \eqref{4.24} imply
\begin{align}\nn
\E_A(x,t)=&\E^{\Om_0(x,t)}(I_3+\Om_1(x,t)A_0+\Om_2(x,t)A_0^2)
\\ \label{4.28}
=&\E^{\Om_0(x,t)}\begin{bmatrix}
1 & \Om_1(x,t) & \Om_2(x,t) 
\\ 0 & 1  & \Om_1(x,t)
\\ 0& 0 & 1 
\end{bmatrix}, 
\\ \nn  
\Om_2(x,t):= & \frac{1}{2} \Om_1(x,t)^2+(\mu_0-c_2)^{-3}x_2+(\mu_0-c_3)^{-3}x_3
\\ \label{4.29} &  
-\I\big(1+(\mu_0-c_2)^{-4}+(\mu_0-c_3)^{-4}\big)t.
\end{align}
\end{Ee}
\begin{Cy}\label{CyN3} Let $k=p=3$, $N=3$, and $ \vk=\mu_0+\ov{\mu_0}>0$.
Define matrices $C_{\Psi}$, $C_{\Phi}$, and $R$ via \eqref{4.26}
and \eqref{4.27}. Choose $\wh C_{\Psi}= \wh C_{\Phi}^*$ and $\cls_0>0$.
Then   relations \eqref{2.6} and \eqref{3.1'}-\eqref{3.3},
where $\E_B(-x, -t)=\E_A(x,t)^*$ and $\E_A$ is given by \eqref{4.28} $($using
functions $\Om_i$, which are introduced in \eqref{4.22},
\eqref{4.23},   \eqref{4.29}$)$,
explicitly define non-singular solutions and potentials of TDS.
\end{Cy}
If we take the $S_3$-node from Example \ref{CE2} and set $\cls_0=0$,
then relations \eqref{4.30} hold. Therefore, taking into account Proposition
\ref{PnRS} we see that the potential $\wt q$ is rational.
Usually, in the study of {\it lumps} it is required that the corresponding potentials
(or solutions) are not only rational but also non-singular.
To choose non-singular $\wt q$, recall that \eqref{4.20}
is valid for $R$ of the form \eqref{4.27}. Hence, in view of Proposition \ref{PnSA}
conditions
\begin{align}\label{4.33}&
\vk=\mu_0+\ov{\mu_0}>0, \quad \rk(\wh C_{\Phi})=n \leq 3
\end{align}
imply that $\wt q$ is non-singular, and our next corollary follows.
\begin{Cy}\label{CyLump} Let $k=p=3$, $N=3$, and $\a=\I$.
Define matrices $C_{\Psi}$, $C_{\Phi}$, and $R$ via \eqref{4.26}
and \eqref{4.27}. Let numbers $\vk, \, \mu_0, n$ and matrix
$\wh C_{\Phi}= \wh C_{\Psi}^*$ satisfy \eqref{4.33}, and set $\cls_0=0$.
Then   relations \eqref{2.6} and \eqref{3.1'}-\eqref{3.3},
where $\E_B(-x, -t)=\E_A(x,t)^*$ and $\E_A$ is given by \eqref{4.28} $($using
functions $\Om_i$, which are introduced in \eqref{4.22},
\eqref{4.23},   \eqref{4.29}$)$,
explicitly define non-singular solutions and potentials of TDS.
Moreover, the potentials $\wt q$ are not only non-singular but also rational.
\end{Cy}

\section{Conclusion} \label{Concl} 
Thus, a version of the binary Darboux transformation is constructed
for TDS equation, where $k\geq 1$. (No binary Darboux transformations
for the TDS equation, where $k>1$, were known before.)
This is an iterated GBDT version. New families of non-singular
and rational potentials and solutions are obtained. Some results are new
for the case that $k=1$ too.

A certain generalization of a colligation introduced by M.S. Liv\v{s}ic and
of the $S$-node introduced by L.A. Sakhnovich, which we call
$S$-multinode, is used in our construction, and could be
useful also in constructions of explicit solutions for other 
multidimensional systems. Another interesting possibility
is application to generalized multidimensional 
nonlinear Schr\"odinger equations 
from \cite{Sab0}.

{\bf Acknowledgement.}
This work was supported by the Austrian Science Fund (FWF) under
Grant  no. Y330.



\end{document}